\documentclass[preprint,12pt]{elsarticle}
\graphicspath{ {./figures/} }
\usepackage{subfig}
\usepackage{xcolor}
\usepackage{amssymb}
\usepackage{amsmath}
\usepackage{hyperref}
\newtheorem{The}{Theorem}[section]
\newtheorem{Def}{Definition}[section]
\newtheorem{Lem}{Lemma}[section]
\newtheorem{Ex}{Example}[section]
\journal{}
\begin{document}
\begin{frontmatter}

\title{Stability analysis of fixed point 
of fractional-order coupled map lattices}

 \author[1]{Sachin Bhalekar}
 \ead{sachin.math@yahoo.co.in, sachinbhalekar@uohyd.ac.in (Corresponding Author)}

 \author[2]{Prashant M. Gade}
 \ead{prashant.m.gade@gmail.com}

 \address [1]{School of Mathematics and Statistics, University of Hyderabad, Hyderabad, 500046 India}
 \address [2]{Department of Physics, RTM Nagpur University, Nagpur}

\begin{abstract}
 We study the stability of synchronized
fixed-point state for linear fractional-order coupled
    map lattice(CML). We observe that the eigenvalues of 
	the connectivity matrix
    determine the stability as for integer-order CML.
    These eigenvalues can be determined exactly in certain
    cases. We find exact bounds in one-dimensional lattice with
    translationally invariant coupling using the theory
    of circulant matrices. This can be extended to any
    finite dimension.  Similar analysis can be carried
out for the synchronized fixed point of nonlinear coupled
fractional maps where eigenvalues of the Jacobian matrix
play the same role. The analysis is generic and demonstrates that
the 
eigenvalues of connectivity matrix
play a pivotal  role in stability analysis of synchronized fixed point
even in coupled fractional maps.
\end{abstract}
\end{frontmatter}

\section{Introduction}
Fractional dynamics extends the dynamical systems to systems with 
memory and studies in fractional order differential equations have exploded
in the recent past. In integer-order systems,  
dynamical systems theory has been enriched by studies in flows as well as
maps. Numerical difficulties are almost absent in the simulation of maps. 
Almost all routes to chaos  observed in flows are observed in maps as
well \cite{ott}. Most chaos control schemes are applicable in the
flows as well as
maps. They appear naturally in scientific contexts where time is discrete.
They have found applications in convecting fluids, lasers, heart cells, 
chemical oscillators, etc \cite{strogatz}.
Circle map has found applications
in several systems described by a damped driven pendulum. Examples 
include Josephson junction in microwave field \cite{bak-josephson},
charge density waves, lasers \cite{laser1,laser2} cardiac arrhythmia
\cite{glass} and even 
air-bubble formation \cite{bubble}. 
Logistic maps have found applications in chemical physics and
population dynamics \cite{chemical,ecology}
These systems have been extended to a spatially
extended version popularly known as coupled map lattice. Coupled map 
lattices have found applications in  diverse
fields such as austenite-martensite
structural transformation, convection and  
crystal growth \cite{shenoy,convection,crystal}. 

Thus it can be useful to investigate coupled fractional maps to
understand the dynamics of spatiotemporal systems in presence of memory.
Studies in fractional maps are  unfrequent compared to fractional
differential equations. Simulation of the fractional differential equation
is computationally cumbersome compared to ordinary differential
equation. It also needs domain expertise in
numerical analysis. Simulating high dimensional system of fractional  
differential equations will need extensive computational 
resources. Though simulation of fractional maps is
more time-consuming than integer-order maps, the computational
resources required are far less than that for fractional
differential equations. 

Systems with power-law memory occur in several physical situations
ranging from electromagnetic waves in dielectric media to adaptation
in biological systems \cite{electromagnetic,bioengineering}.
In this work, we study coupled fractional maps and investigate 
a very basic problem of existence and stability of fixed-point
solution and state conditions for a synchronized fixed-point solution. 
In certain important cases, such as coupled map lattice 
in finite dimension, explicit bounds can be derived for stability.

\section{Preliminaries} \label{prel}
In this section, we present some basic definitions and results.
Let $ h > 0 ,\; a \in \mathbb{R}$, $ (h\mathbb{N})_a = \{ a, a+h, a+2h, \ldots\} $ and $\mathbb{N}_a=\{a, a+1, a+2, \ldots \}$.
\begin{Def}(see \cite{ferreira2011fractional, bastos2011discrete, mozyrska2015transform}).
	For a function $x : (h\mathbb{N})_a \rightarrow  \mathbb{R}$, the forward h-difference operator if defined as 
	$$
	(\Delta_h x)(t)=\frac{x(t+h)- x(t)}{h},$$
	where t	$ \in (h\mathbb{N})_a $.
\end{Def}
	Throughout this paper, we take $a = 0$ and $h = 1$.
\begin{Def}\cite{mozyrska2015transform}
	For a function  $x : \mathbb{N_\circ} \rightarrow  \mathbb{R}$ the fractional sum of order $\alpha>0$ is given by
	\begin{equation}
	(\Delta^{-\alpha}x)(t) 
	= \frac{1}{\Gamma(\alpha)}\sum_{s=0}^{n}\frac{\Gamma(\alpha+n-s)}{\Gamma(n-s+1)} x(s),
	\end{equation}
	where,	$t=\alpha+n, \; n \in \mathbb{N_\circ}$.
\end{Def}

\begin{Def}\cite{mozyrska2015transform,fulai2011existence}
	Let $\mu>0$ and $m-1<\mu<m$, where $m\in \mathbb{N}$, $m=\lceil \mu \rceil$. The $\mu$th fractional Caputo like difference is defined as
	\begin{equation}
	\Delta^\mu x(t)= \Delta^{-(m-\mu)}\left(\Delta^m x(t)\right),
	\end{equation}
	where $t\in \mathbb{N}_{m-\mu}$ and 
	\begin{equation}
	\Delta^m x(t)=\sum_{k=0}^{m}\binom{m}{k}(-1)^{m-k}x(t+k).
	\end{equation}
\end{Def}

	\begin{Def}	 \cite{mozyrska2015transform} 
		The Z-transform of a sequence $ \{y(n)\}_{n=0}^\infty $ is a complex function given by
	$Y(z)=Z[y](z)=\sum_{k=0}^{\infty} y(k) z^{-k}$
	where $z \in \mathbb{C}$ is a complex number for which the series converges absolutely.
\end{Def}

\begin{Def}\cite{mozyrska2015transform} 
Let $\tilde{\phi}_\alpha(n)$ be a family of binomial functions defined on $\mathbb{Z}$, parametrized by $\alpha$ defined by
\begin{eqnarray}
\tilde{\phi}_\alpha(n) &=& \frac{\Gamma(n+\alpha-1)}{\Gamma(\alpha) \Gamma(n)}\nonumber\\
&=& \left(
\begin{array}{c}
n+\alpha-1\\
n\\
\end{array}
\right)
=(-1)^n
\left(
\begin{array}{c}
-\alpha\\
n
\end{array}
\right).
\end{eqnarray}
Then
\begin{equation*}
Z(\tilde{\phi}_{\alpha}(t))=\frac{1}{(1-z^{-1})^{\alpha}}, \quad |z|>1.
\end{equation*}
\end{Def}

\begin{Def}  \cite{mozyrska2015transform} 
The convolution $\phi*x$ of the functions $\phi$ and $x$ defined on $\mathbb{N}$ is defined as
\begin{equation*}
\left(\phi*x\right)(n)=\sum_{s=0}^{n}\phi(n-s)x(s)=\sum_{s=0}^{n}\phi(s)x(n-s).
\end{equation*}
Then the Z-transform of this convolution is 
\begin{equation}
Z\left(\phi*x\right)(n)=\left(Z\left(\phi\right)(n)\right) \left(Z\left(x\right)(n)\right).
\end{equation}
\end{Def}	

\begin{Lem}\cite{fulai2011existence}
	The discrete function $x(t)$ is solution of an initial value problem
	\begin{eqnarray}
		\Delta^\alpha x(t)&=& f(x(t+\alpha-1)), \quad t\in N_{1-\alpha},\quad 0<\alpha<1,\nonumber\\
		x(0)&=&x_0
	\end{eqnarray}
	if and only if $x(t)$ is solution of following fractional discrete dynamical system
	\begin{eqnarray}
	x(t)&=&x_0+\sum_{s=1-\alpha}^{t-\alpha}\frac{\Gamma(t-s)}{\Gamma(\alpha) \Gamma(t-s-\alpha+1)}f\left(x(s+\alpha-1)\right)\nonumber\\
	&=&x_0+\sum_{j=0}^{t-1}\frac{\Gamma(t-j+\alpha-1)}{\Gamma(\alpha) \Gamma(t-j)}f\left(x(j)\right).
	\end{eqnarray}
\end{Lem}

\section{Fractional order coupled map lattices: Linear systems}\label{lin}
Consider the linear coupled map lattice of fractional order $\alpha\in(0, 1)$
\begin{equation}
x_{t+1}(k)=x_0(k)+
\sum_{j=0}^{t}
\sum_{m=1}^N\frac{\Gamma(t-j+\alpha)}{\Gamma(\alpha) \Gamma(t-j+1)}
\left(A_{km} x_j(m)-x_j(k)\right), \label{linCML}
\end{equation}
where $x_t(k)$ is the variable at time $t$ associated with the $k$-th lattice 
point, $k=1,2,\cdots,N$, $x_t(0)=x_t(N)$ and $x_t(N+1)=x_t(1)$
and $A=(A_{km})$ is $N\times N$ connectivity matrix.\\
If we write $X_t=\left(x_t(1), x_t(2), \cdots, x_t(N)\right)$, 
a column vector in $\mathbb{R}^N$ 
then the system (\ref{linCML}) is equivalent to 
\begin{eqnarray}
X_{t+1}&=&X_0 + \sum_{j=0}^{t}\frac{\Gamma(t-j+\alpha)}{\Gamma(\alpha) \Gamma(t-j+1)} \left(A-I\right)X_j,\nonumber\\
&=&X_0 + \left(A-I\right) \left(\tilde{\phi}_{\alpha}(t)*X_t\right),
 \label{linCMLSys}
\end{eqnarray}
where $I$ is $N\times N$ identity matrix. Applying Z-transform and using the properties given in Section \ref{prel}, we get
\begin{equation*}
z \bar{X}(z) -z X_0=\frac{1}{1-z^{-1}}X_0+\frac{1}{(1-z^{-1})^{\alpha}}\bar{X}(z)\left(A-I\right), \quad |z|>1,
\end{equation*}
where $\bar{X}(z)$ is the Z-transform of $X_t$. Therefore, the 
characteristic equation of system  (\ref{linCMLSys}) is given as
\begin{equation}
det\left(z(1-z^{-1})^{\alpha}I-(A-I)\right)=0. \label{cheq}
\end{equation}
Motivated from 
\cite{stanislawski2013stability,vcermak2015explicit}, 
we propose the following stability theorem.
\begin{The} \label{stabth0}
The zero solution of the system (\ref{linCML}) or (\ref{linCMLSys}) is asymptotically stable if and only if all the roots of the characteristic equation (\ref{cheq}) satisfy $|z|<1$.	
\end{The}

\subsection{Stable Region}
\begin{figure}%
	\centering
	\includegraphics[width=10cm]{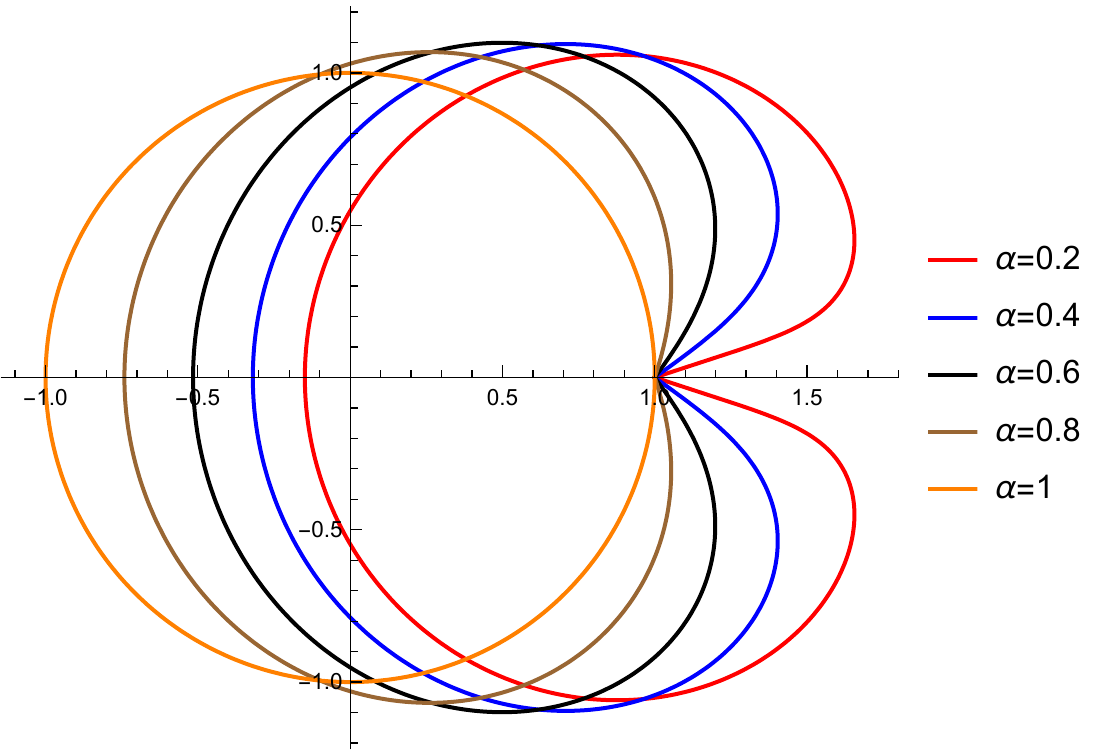} 
	\caption{The stability region of fractional order map}%
	\label{fig1}%
\end{figure}

At the boundary of stable region, the root $z$ of characteristic equation  (\ref{cheq}) should satisfy $|z|=1$. Therefore, we obtain the parametric boundary curve $\beta(t)$ of stable region by substituting $z=e^{\iota t}$, $0\leq t \leq 2\pi$ in the (\ref{cheq}) as
\begin{equation}
\beta(t)=\left(2^\alpha \left(\sin(t/2)\right)^\alpha \cos\left(\alpha \frac{\pi}{2}+t(1-\alpha/2)\right)+1, 2^\alpha \left(\sin(t/2)\right)^\alpha \sin\left(\alpha \frac{\pi}{2}+t(1-\alpha/2)\right) \right). \label{bdry}
\end{equation}
The boundary curves $\beta(t)$ for different values of $\alpha \in (0, 1]$ are sketched in Figure \ref{fig1}.
If the eigenvalues of $A$ are complex, we need to consider if the 
given eigenvalue is in the stable region defined
by the cardioid given above and the solution is stable only if
all eigenvalues lie in the stable region.\\

We have following result \cite{stanislawski2013stability,gade2021fractional,buslowicz2013necessary,vcermak2015explicit}.
\begin{The} \label{stabth1}
If all the eigenvalues of matrix $A$ lie inside the region bounded by the curve $\beta(t)$, $0\leq t \leq 2\pi$ defined by (\ref{bdry}) then the system (\ref{linCML}) is asymptotically stable.
\end{The}

Thus the stability of synchronized fixed point $x_t(k)=0$ as $t\rightarrow
\infty$ $\forall i$ depends only on eigenvalues of connectivity matrix $A$.
Let us consider a particular case of coupled map lattice on 
one-dimensional lattice with translationally invariant
coupling and periodic boundary conditions. The matrix $A$ such that 
$A_{ii}=a_1$, $A_{i,i+1}=a_2$
and $A_{i,i-1}=a_0$ is given by.
\begin{equation*}
A=\begin{pmatrix}
a_1 & a_2 & 0 & 0 & \cdots & 0 & 0 & a_0\\
a_0 & a_1 & a_2 & 0 & \cdots & 0 & 0 & 0\\
0 & a_0 & a_1 & a_2 & \cdots & 0 & 0 & 0\\
\vdots & \vdots & \vdots & \vdots &  & \vdots & \vdots & \vdots\\
0 & 0 & 0 & 0 & \cdots &a_0 & a_1 & a_2\\
a_2 & 0 & 0 & 0 & \cdots & 0 & a_0 & a_1\\
\end{pmatrix}
\end{equation*}
For the special case $N=2$, we define
\begin{equation*}
A=\begin{pmatrix}
a_1 & a_0+a_2\\
a_0+a_2 & a_1 \\
\end{pmatrix}
\end{equation*}
and for $N=1$
\begin{equation*}
A=\begin{pmatrix}
a_1 + a_0+a_2
\end{pmatrix}
\end{equation*}
using periodic boundary conditions.

$A$ is a circulant matrix with eigenvalues
$\lambda_l=a_1+a_2 \omega^l+ a_0 \omega^{-l}$ where
$\omega=\exp(\iota {\frac{2\pi}{N}})$
 is primitive $N$th root of unity \cite{amritkar1991stability}. 
For symmetric case, where
$a_2=a_0$, we get $\lambda_l=a_1+2a_2 \cos(\theta_l)$,
where $\theta_l={\frac{2 \pi l}{N}}$
 for $0\le l\le N-1$.
We note that $\lambda_l=\lambda_{N-l}$ in this case.
For
case $a_2=-a_0$, we obtain $\lambda_l=a_1+\iota 2a_2 \sin(\theta_l)$
where $\theta_l={\frac{2 \pi l}{N}}$.
If $N=4K$, $\lambda_K=a_1+\iota 2a_2$ and $\lambda_{3K}=a1-\iota 2a2$.
These are limiting cases in this system.
Coupled map lattice in one dimension is a widely explored system and
we will study the above cases in further detail in this section.

First we consider the bounds on real part of the eigenvalue.\\
\textbf{Note:} The stable region of the real eigenvalue 
$\lambda$ is $1-2^{\alpha}<\lambda<1$.\\

\begin{Ex}
	Consider  $N=3$ and $\alpha=0.4$. The parameter values  $a_0=0.2$,  $a_1=-0.5$ and  $a_2=0.1$ produce the eigenvalues $-0.2$ and $-0.65 \pm 0.0866025 \iota$. Since the eigenvalues $-0.65 \pm 0.0866025 \iota$ are outside the stable region, we get the unstable solutions as shown in Fig. \ref{fig2}.
	\begin{figure}%
		\centering
		\includegraphics[width=10cm]{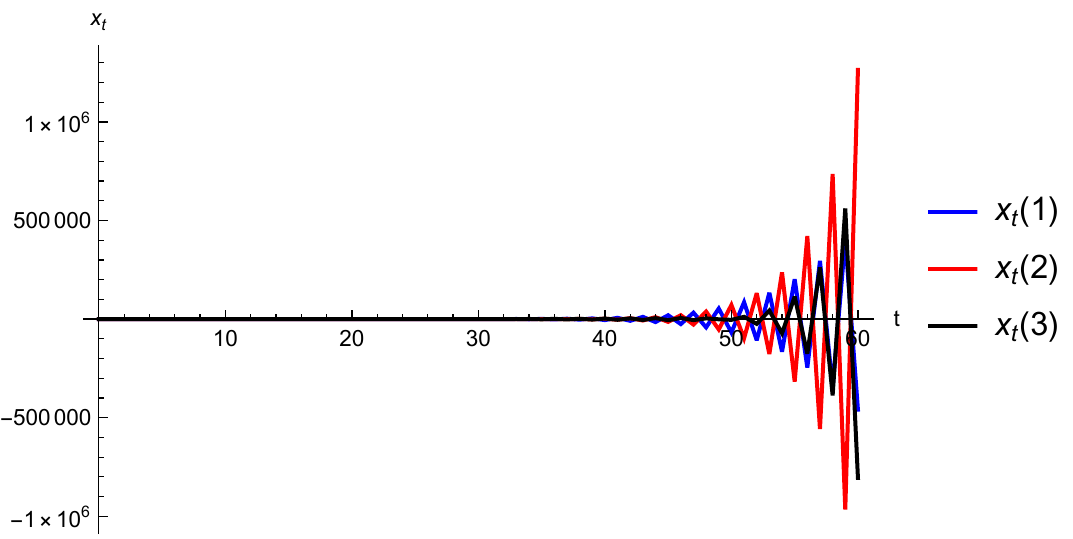} 
		\caption{Unstable solution for $N=3$, $\alpha=0.4$,  $a_0=0.2$,  $a_1=-0.5$ and  $a_2=0.1$}%
		\label{fig2}%
	\end{figure}
On the other hand, if we set  $N=3$, $\alpha=0.8$,  $a_0=0.2$,  $a_1=-0.3$ and  $a_2=0.1$ then all the eigenvalues viz. $0$ and $-0.45 \pm 0.0866 \iota$ lie inside the stable region and we get the stable solutions (cf. Fig. \ref{fig3}).
\begin{figure}%
	\centering
	\includegraphics[width=10cm]{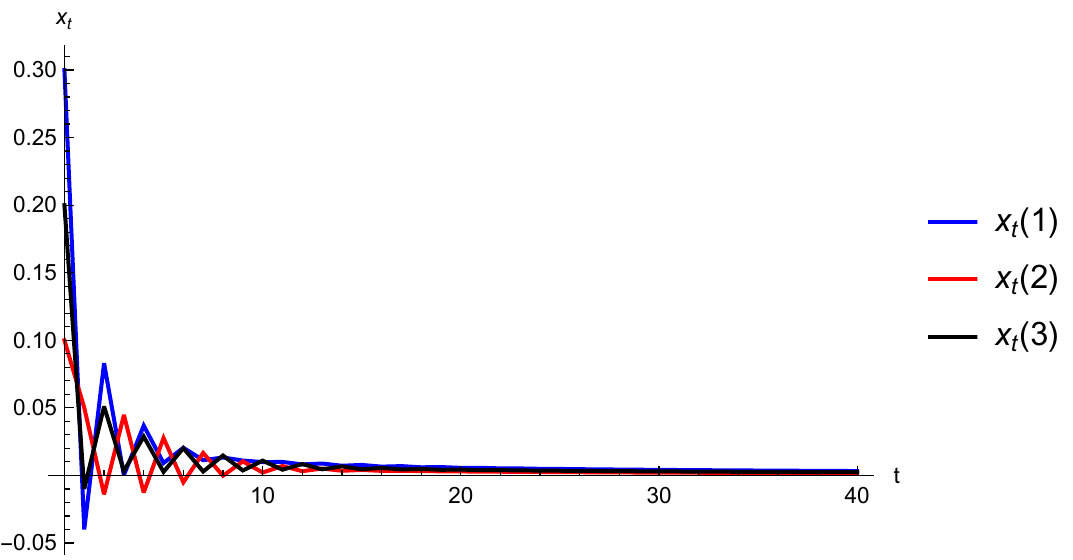} 
	\caption{Stable solution for $N=3$, $\alpha=0.8$, $a_0=0.2$,  $a_1=-0.3$ and  $a_2=0.1$}%
	\label{fig3}%
\end{figure}
\end{Ex}

In the following two subsections, we consider two important particular cases and discuss the stability.
\subsubsection{Symmetric Case}
\begin{figure}%
	\centering
	\includegraphics[width=10cm]{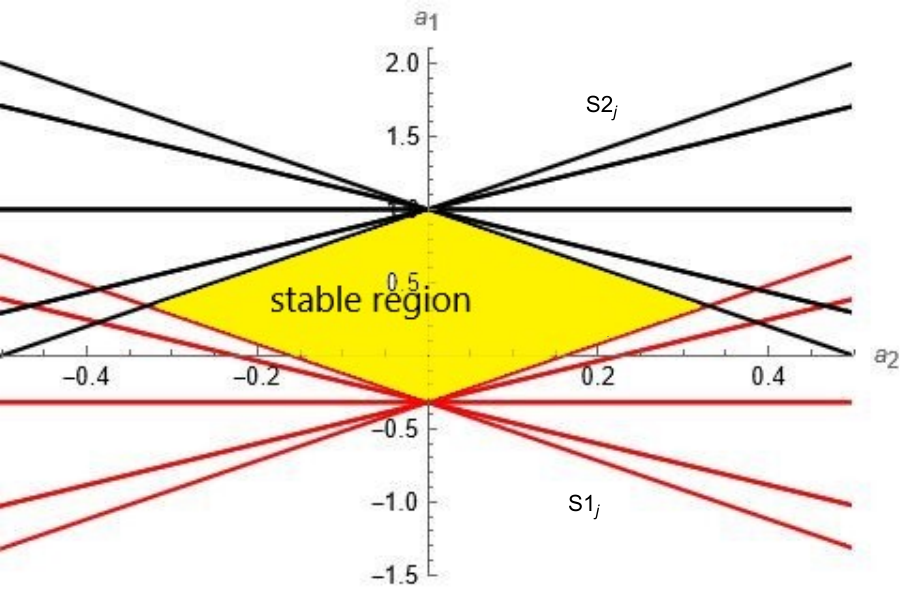} 
	\caption{Stable region of symmetric case}%
	\label{fig4}%
\end{figure}

In this section, we assume that the coefficient matrix $A$ in (\ref{linCMLSys}) is symmetric i.e. $a_0=a_2$.
 \begin{The}\label{sym}
	The stable region of the system (\ref{linCMLSys}) with $a_0=a_2$ is bounded by the quadrilateral with vertices $Q_1=(0,1)$, $Q_2=\left(-2^{\alpha-2}, 2^{\alpha-1}+1-2^{\alpha}\right)$, $Q_3=(0,1-2^\alpha)$ and $Q_4=\left(2^{\alpha-2}, 1-2^{\alpha-1}\right)$ for even values of lattice points $N$ and $Q_1=(0,1)$, $Q_2'=\left(-\frac{2^{\alpha-1}}{1+\cos(\pi/N)}, \frac{2^{\alpha}}{1+\cos(\pi/N)}+1-2^{\alpha}\right)$, $Q_3=(0,1-2^\alpha)$ and $Q_4'=\left(\frac{2^{\alpha-1}}{1+\cos(\pi/N)}, -\frac{2^{\alpha}}{1+\cos(\pi/N)}+1\right)$ for odd $N$ in the   $a_2a_1$-plane.
\end{The}
{\textbf{Proof:}}
 The eigenvalues of $A$ in the symmetric case are 
\begin{equation}
\lambda_j=a_1+2 a_2 \cos\left(\frac{2\pi j}{N}\right), \quad j=0,1,\cdots,N-1.
\end{equation}
Note that $\cos\left(\frac{2\pi j}{N}\right)=\cos\left(\frac{2\pi (N-j)}{N}\right)$. Therefore, to obtain the distinct values we take $j=0,1,\cdots,[N/2]$, where $[r]$ is an integer-part of the real number $r$.
Since, all these eigenvalues are real, the stable region in the $a_2a_1$-plane is an intersection of the regions
\begin{equation}
1-2^{\alpha}<a_1+2 a_2 \cos\left(\frac{2\pi j}{N}\right)<1, \quad j=0,1,\cdots,[N/2]. \label{regns}
\end{equation}
The boundaries of these regions are straight lines defined by following two sets
\begin{eqnarray}
S1_j: \quad a_1 &=& -2\cos\left(\frac{2\pi j}{N}\right) a_2 + (1-2^{\alpha}), \quad \text{and} \label{line1} \\
S2_j: \quad a_1 &=& -2\cos\left(\frac{2\pi j}{N}\right) a_2 + 1, \label{line2} 
\end{eqnarray}
where $j=0,1,\cdots,[N/2]$.\\
 Note that
\begin{equation}
[N/2]=
\begin{cases}
N/2, \quad \text{if $N$ is even}\\
(N-1)/2, \quad \text{if $N$ is odd.}
\end{cases}
\end{equation}
The stable region of (\ref{linCMLSys}) will be bounded by the straight lines $S1_j$ and $S2_j$ which are close to origin in the $a_2a_1$-plane, as shown in Figure \ref{fig4}.
 The lines in the set $S1_j$ intersect each other at $(0,1-2^\alpha)$ whereas those in the set $S2_j$ have intersection at  $(0,1)$ in the $a_2a_1$-plane. 
 Further, the innermost line $a_1 = -2 a_2 + (1-2^{\alpha})$ in the set (\ref{line1}) with $j=0$ intersects the innermost line $a_1 = -2\cos\left(\frac{2\pi [N/2]}{N}\right) a_2 + 1$ in the set (\ref{line2}) with $j=[N/2]$ in the $a_2a_1$-plane at the point $\left(-2^{\alpha-2}, 2^{\alpha-1}+1-2^{\alpha}\right)$ when $N$ is even and at $\left(-\frac{2^{\alpha-1}}{1+\cos(\pi/N)}, \frac{2^{\alpha}}{1+\cos(\pi/N)}+1-2^{\alpha}\right)$ when $N$ is odd. Secondly, the intersection between the innermost lines $a_1 = -2\cos\left(\frac{2\pi [N/2]}{N}\right) a_2 + (1-2^{\alpha})$ in the set (\ref{line1}) with  $j=[N/2]$ and $a_1 = -2 a_2 + 1$ in the set (\ref{line2}) with $j=0$ is the point $\left(2^{\alpha-2}, 1-2^{\alpha-1}\right)$ when $N$ is even and $\left(\frac{2^{\alpha-1}}{1+\cos(\pi/N)}, -\frac{2^{\alpha}}{1+\cos(\pi/N)}+1\right)$ when $N$ is odd.\\
 Thus, the stable region which is an intersection of all the regions (\ref{regns}) is bounded by the quadrilateral with vertices described in the statement of this theorem.  This proves the result.\\

We note that stable region does not change for even $N$. 
Two extreme values for $a_1+2a_2\cos(\theta_l)$ are given by
$\lambda_0=a_1+2a_2$ and $\lambda_{N/2}=a_1-2a_2$.
For odd $N$, one of the limits $\lambda_0=a_1+2a_2$ is
still realized. Other limit is slightly bigger by a leading to a slightly
higher stability range and the it is approached as $1/N^2$ for
large $N$.  
Thus, in the 
thermodynamic limit $N\rightarrow \infty$, 
stability region for $N\rightarrow \infty$
coincides with stability region for $N=2$. 
Thus  
the stability of extreme eigenvalues in the thermodynamic
limit determine the stability region.

\begin{Ex}
Consider the symmetric system  (\ref{linCMLSys}) with even number $N=8$ of lattice points and $\alpha=0.2$. The stable region using Theorem \ref{sym} is sketched in Figure \ref{fig5}. 
\begin{figure}%
	\centering
	\includegraphics[width=10cm]{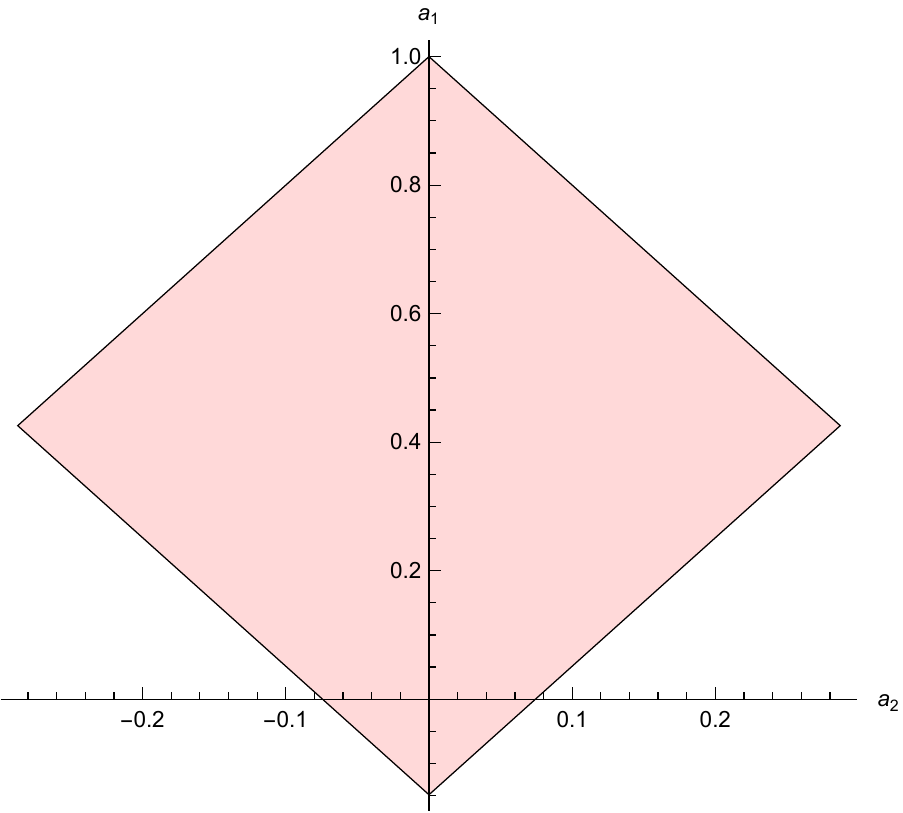} 
	\caption{Stable region of symmetric  system  (\ref{linCMLSys}) with $N=8$ and $\alpha=0.2$}%
	\label{fig5}%
\end{figure}
We verified that the solutions starting in a neighborhood of origin converge to origin if we take $(a_2, a_1)$ in the stable region. Figure \ref{fig6} shows the converging trajectories for the parameter values $a_1=0.1$, and $a_2=-0.05$. The unstable solution is sketched in Figure \ref{fig7} with $a_1=-0.02$, and $a_2=0.1$ which are outside the stable region.
\begin{figure}%
	\centering
	\includegraphics[width=10cm]{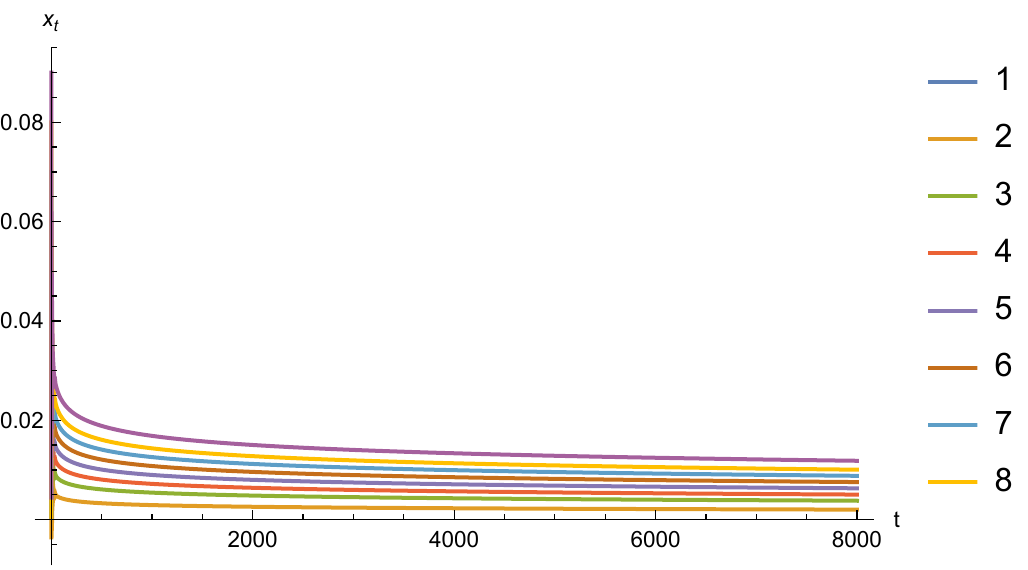} 
	\caption{Stable solution of symmetric  system (\ref{linCMLSys}) with $N=8$,  $\alpha=0.2$, $a_1=0.1$, and $a_2=-0.05$}%
	\label{fig6}%
\end{figure} 
\begin{figure}%
	\centering
	\includegraphics[width=10cm]{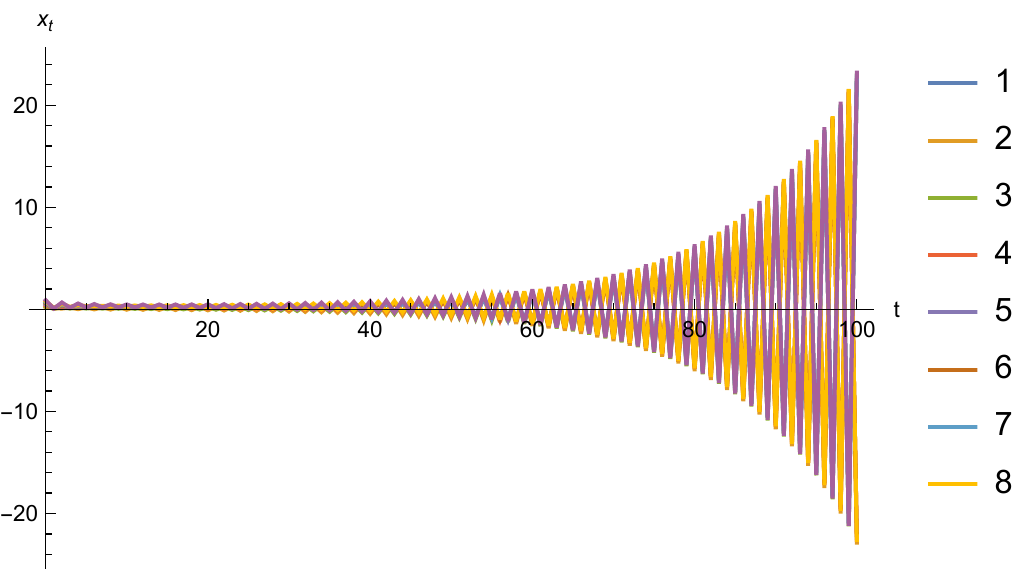} 
	\caption{Unstable solution of symmetric  system (\ref{linCMLSys}) with $N=8$,  $\alpha=0.2$, $a_1=-0.02$, and $a_2=0.1$}%
	\label{fig7}%
\end{figure} 
\end{Ex}
 
 \begin{Ex}
Let us consider the symmetric system  (\ref{linCMLSys}) with odd number $N=9$ of lattice points. The stable region in this case with $\alpha=0.5$ is shown in Figure \ref{fig8}. The stable solution for the parameter values  $a_1=0.6$, and $a_2=-0.1$is shown in Figure \ref{fig9} whereas the unstable solution for $a_1=0.2$, and $a_2=0.6$ is in Figure \ref{fig10}.
\begin{figure}%
	\centering
	\includegraphics[width=10cm]{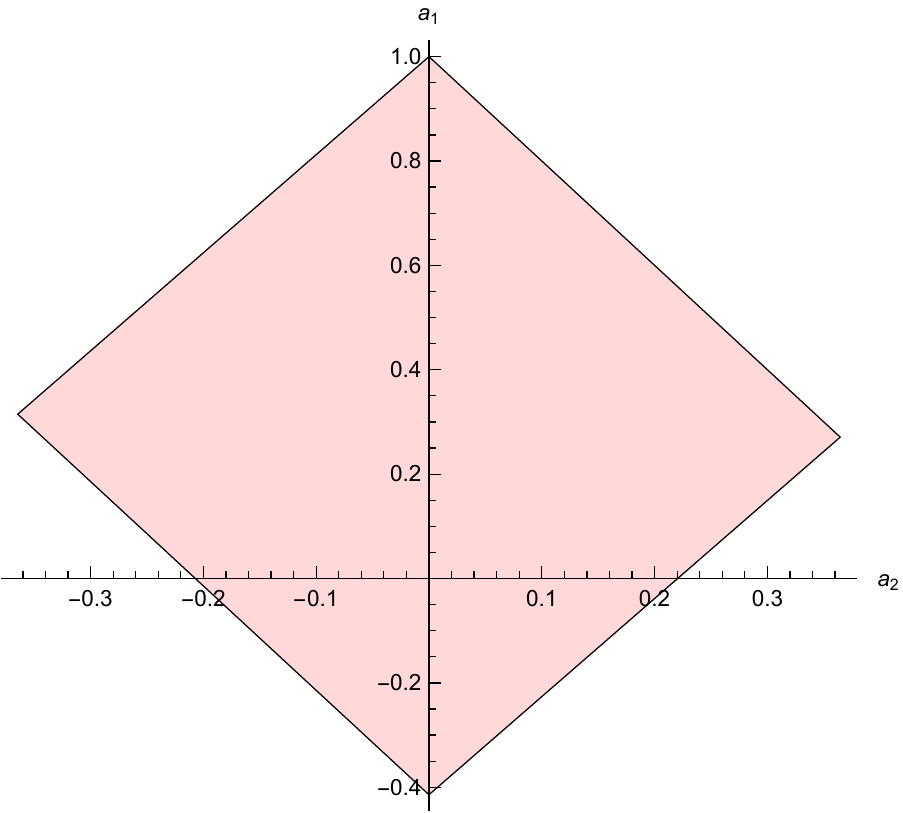} 
	\caption{Stable region of symmetric  system  (\ref{linCMLSys}) with $N=9$ and $\alpha=0.5$}%
	\label{fig8}%
\end{figure}
\begin{figure}%
	\centering
	\includegraphics[width=10cm]{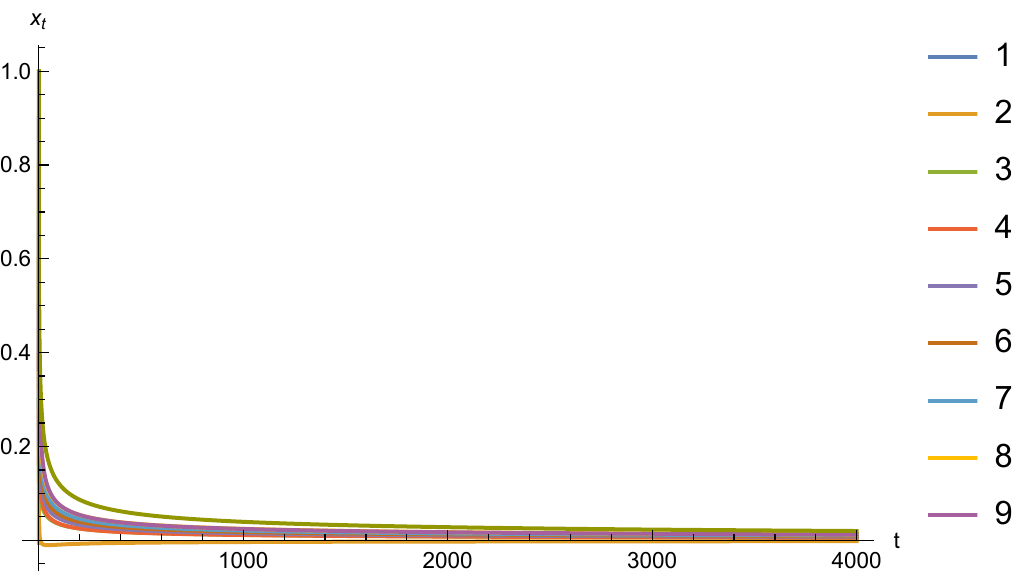} 
	\caption{Stable solution of symmetric  system (\ref{linCMLSys}) with $N=9$,  $\alpha=0.5$, $a_1=0.6$, and $a_2=-0.1$}%
	\label{fig9}%
\end{figure} 
\begin{figure}%
	\centering
	\includegraphics[width=10cm]{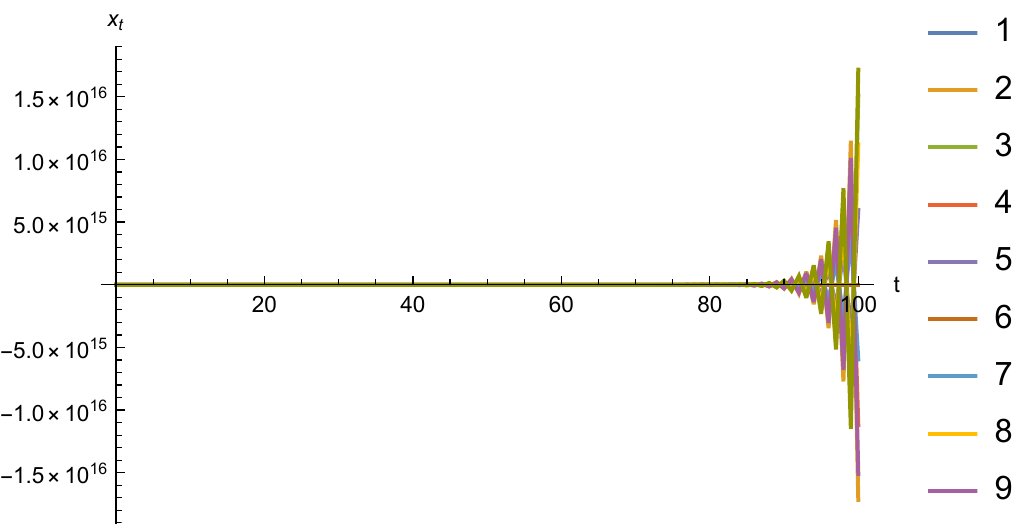} 
	\caption{Unstable solution of symmetric  system (\ref{linCMLSys}) with $N=9$,  $\alpha=0.5$, $a_1=0.2$, and $a_2=0.6$}%
	\label{fig10}%
\end{figure}
 \end{Ex}

\subsubsection{Asymmetric Case}
Now, we consider the asymmetric system (\ref{linCMLSys}) with  $a_0=-a_2$. We define the cardioids 
\begin{equation}
\gamma_j=\left(Re\left[e^{\iota t} \left(1-e^{-\iota t}\right)^{\alpha} \right]+1, \frac{1}{2\sin(2\pi j/N)}Im\left[e^{\iota t} \left(1-e^{-\iota t}\right)^{\alpha} \right]\right), \quad j=1,2,\cdots,[N/2] \label{gam}
\end{equation}
in the $a_1a_2$-plane, provided $\sin(2\pi j/N)\neq 0$.\\
 The stability result in this case is discussed below. Note that $\lceil r \rceil$ is the ceiling function of real number $r$.

\begin{figure}%
	\centering
	\includegraphics[scale=1]{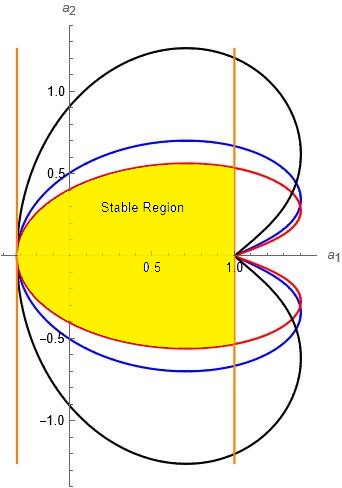} 
	\caption{Stable region of asymmetric system with $a_0=-a_2$}%
	\label{fig11}%
\end{figure}

\begin{The}\label{asymStabr}
Consider the system (\ref{linCMLSys}) with  $a_0=-a_2$. We have following stability results:
\begin{itemize}
	\item If $N=1$ or $N=2$ then the stable region is $1-2^\alpha<a_1<1$.
	\item If $N\geq 3$ is an odd number then the stable region is bounded by the line $a_1=1$ and the cardioid $\gamma_{\lceil \frac{N-1}{4}\rceil}$ in the $a_1a_2$-plane.
	\item If $N\geq 4$ is an even number then the stable region is bounded by the line $a_1=1$ and the cardioid $\gamma_{\left[ \frac{N}{4}\right]}$ in the $a_1a_2$-plane.
\end{itemize}
\end{The}
{\textbf{Proof:} }
Suppose that $a_0=-a_2$ in (\ref{linCMLSys}) . Therefore, the eigenvalues of the coefficient matrix $A$ are
\begin{equation}
\lambda_j=a_1+\iota 2 a_2 \sin\left(\frac{2\pi j}{N}\right), \quad j=0,1,\cdots,N-1.
\end{equation}
It is observed that for $j=[N/2]+1, [N/2]+2,\cdots,N-1$, the values $\lambda_j$ are complex conjugates of those for $j=0,1,\cdots,[N/2]$.
Therefore, the stable region corresponding to $\lambda_j$ is given by the cardioid $\gamma_j$ defined in (\ref{gam}).\\
Since, $\lambda_0=a_1\in \mathbb{R}$, one of the stability conditions is
\begin{equation}
1-2^\alpha<a_1<1. \label{stab0}
\end{equation}
Further, if $N=1$ or $N=2$ then $a_1$ is the only eigenvalue of matrix $A$. Therefore, the stability condition is given by (\ref{stab0}).\\
For $N\geq 3$, the stable region is an intersection of the cardioids $\gamma_j$ and the region (\ref{stab0}). It is observed that, this region is bounded by the ``innermost" cardioid and the line $a_1=1$ in the $a_1a_2$-plane, as shown in the Figure \ref{fig11}.
Now, we have to find the $j$ for which the cardioid $\gamma_j$ is innermost.\\
The innermost cardioid is generated by $\gamma_j$ for which the value $\sin(2\pi j/N)$ is maximum. Further, the value $\sin(2\pi j/N)$ is maximum for the number $2\pi j/N$ which is closest to $\pi/2$ i.e. if the value $\vert\frac{2\pi j}{N}-\frac{\pi}{2}\vert = \frac{\pi}{2N}\vert 4j-N\vert$ is minimum.\\
Thus, our problem is reduced to find minimum of the set
\begin{equation}
S=\left\{\vert 4j-N\vert: j=1,2,\cdots,[N/2]\right\}.
\end{equation}
If $N$ is even number, then the minimum of $S$ occurs at $j=[N/4]$. On the other hand, if $N$ is an odd number, then the minimum of $S$ occurs at $j=\lceil \frac{N-1}{4}\rceil$.\\
The result is proved.

We note that if the number of maps is multiple of 4, say $N=4K$, 
$\lambda_K=a_1+\iota 2 a_2$ and $\lambda_{N-K}=a_1-\iota 2 a_2$. 
Also, 
$\lambda_0=a_1$ for any $N$. The cardioid $\gamma_K$ defined
above reduces to cardioid for given value of $\alpha$ for $j=0$
where real part is given by $a_1$ and imaginary part is $2a_2$.
For $j=0$, $\gamma_j$ is strip between $1-2^{\alpha}\le a_0 \le 1$ with no condition on $a_2$.
The stability region is given
by intersection of this strip with the cardioid $\gamma_K$ for given $\alpha$ for
$N=4K$. 
If $N$ is not an exact multiple of 4, the stability region 
is slightly bigger and as expected it shrinks to stability
region for $N=4$ in the thermodynamic limit.

\begin{Ex}
We consider the system  (\ref{linCMLSys}) with  $a_0=-a_2$, $N=6$ and $\alpha=0.3$. Here, $N$ is even and $[N/4]=1$. According to Theorem \ref{asymStabr}, the stable region is bounded by the cardioid $\gamma_1$ and the line $a_1=1$ in the $a_1a_2$-plane, as shown in Figure \ref{fig12}.
\begin{figure}%
	\centering
	\includegraphics[width=10cm]{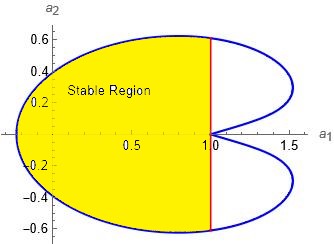} 
	\caption{Stable region of system (\ref{linCMLSys}) with $a_0=-a_2$, $N=6$ and $\alpha=0.3$}%
	\label{fig12}%
\end{figure}
The point $a_1=-0.3, a_2=0.5$ is outside the stable region and therefore we get unstable solution (cf. Figure \ref{fig13}). On the other hand, we get the stable solution (cf. Figure \ref{fig14}) for the parameter values $a_1=-0.1, a_2=-0.22$.
\begin{figure}%
	\centering
	\includegraphics[width=10cm]{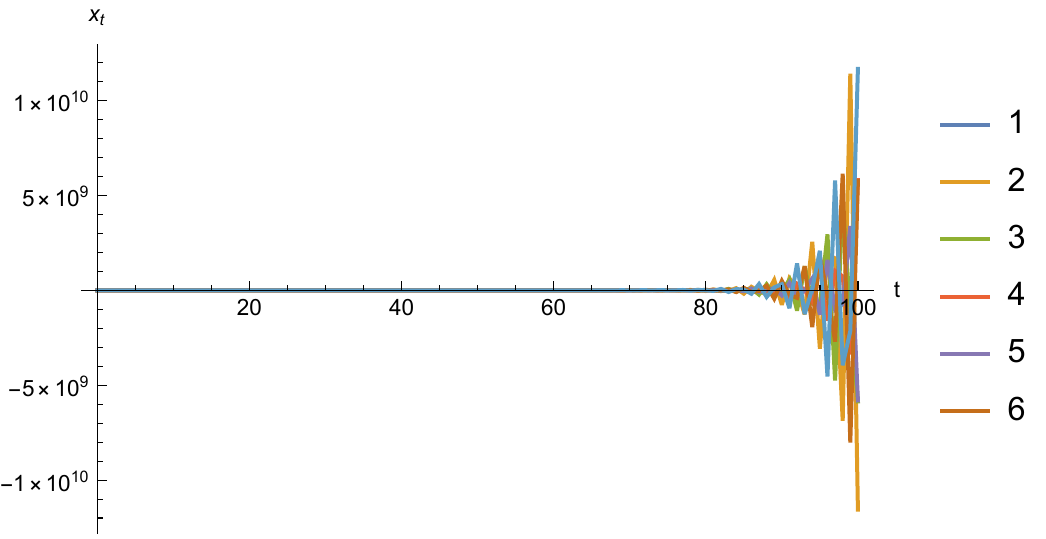} 
	\caption{Unstable solution of (\ref{linCMLSys}) with $a_1=-0.3, a_2=0.5$, $a_0=-a_2$, $N=6$ and $\alpha=0.3$}%
	\label{fig13}%
\end{figure}
\begin{figure}%
	\centering
	\includegraphics[width=10cm]{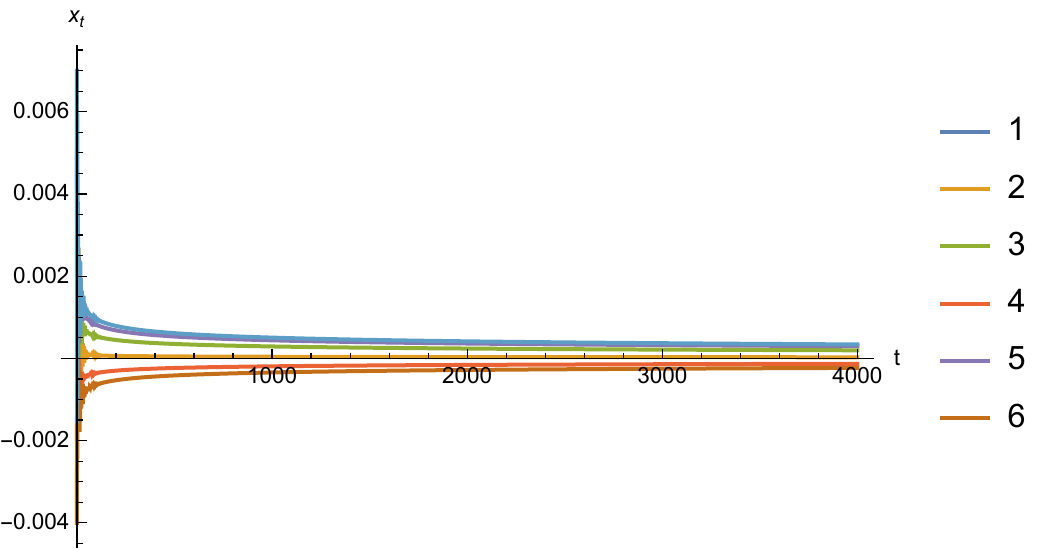} 
	\caption{Stable solution of system (\ref{linCMLSys}) with $a_1=-0.1, a_2=-0.22$, $a_0=-a_2$, $N=6$ and $\alpha=0.3$}%
	\label{fig14}%
\end{figure}
\end{Ex}
We also verified the Theorem \ref{asymStabr} for odd values of $N$ but not presented the example for brevity.

\subsubsection{Thermodynamic limit}

The coupled map lattice (\ref{linCMLSys}) 
in the thermodynamic limit $N\longrightarrow \infty$ 
is an interesting system studied in the literature 
\cite{pikovsky1994collective, carretero1999thermodynamic, 
parekh1998synchronization}. This limit may gives rise 
to some important phenomena such as rescaling of the 
Lyapunov spectrum \cite{carretero1999scaling}. 
The physical interpretation of this limit 
\cite{kaneko1989towards,ruelle1982large} is that a 
coupled map lattice with a very large number of lattice 
points may be identified as a chain of relatively 
small-sized independently evolving subsystems. 
\par As $N\longrightarrow \infty$, $\cos(\pi/N)\longrightarrow 1$. 
Therefore, the stable region of the symmetric system 
(\ref{linCMLSys}) in the $a_1a_2$-plane  according to Theorem \ref{sym} is 
bounded by the quadrilateral $Q_1Q_2Q_3Q_4$ in the thermodynamic limit.
\par As $N\longrightarrow \infty$, the interval $[0, 2\pi]$ 
will contain an infinitely many values of the form $2\pi j/N$, 
$j=1,2,\cdots,N-1$. Therefore, the maximum of $\sin(2\pi j/N)$ will 
approach to $1$ as $N\longrightarrow \infty$. Therefore, the 
stable region in the $a_1a_2$-plane of the asymmetric 
system (\ref{linCMLSys}) with $a_0=-a_2$ 
according to Theorem \ref{asymStabr} is bounded by the 
line $a_1=1$ and the cardioid
\begin{equation*}
\gamma_\infty=\left(Re\left[e^{\iota t} \left(1-e^{-\iota t}\right)^{\alpha} \right]+1, \frac{1}{2}Im\left[e^{\iota t} \left(1-e^{-\iota t}\right)^{\alpha} \right]\right)
\end{equation*}
in the thermodynamic limit.

Though we have studied 1-dimensional case in detail, the formulation is
very generic and can be extended to any case where the eigenvalues of 
the connectivity matrix can be computed analytically. Consider
a 2-dimensional case with $NM$ maps with couplings
$A_{(i,j),(i\pm 1,j)}=a_0$
$A_{(i,j),(i,j\pm 1)}=a_2$
$A_{(i,j),(i,j)}=a_1$. This is  a block-circulant matrix
with circulant blocks and the eigenvalues are given
by $\lambda_{k1,k2}=a_1+2a_0\cos(\theta_{k1})+2a_2\cos(\theta_{k2})$
where $\theta_{k1}={\frac{2\pi k_1}{N}}$
and $\theta_{k2}={\frac{2\pi k_2}{M}}$. The indices $k_1$ and
$k_2$ run from 0 to $N-1$ and 0 to $M-1$ respectively \cite{gade1993spatially}.
The bounds are given by $a_1+2a_0+2a_2$ and $a_1-2a_0-2a_2$ in the thermodynamic
limit (assuming all off-diagonal couplings positive) and the stability region is given by quadrilateral where both
these bounds are in the range $[-2^{\alpha}+1,1]$.
Thus the formulation allows us to analytically
find the stability of a coupled map lattice
with any connectivity matrix if the eigenvalues can be determined analytically.
If we couple each site to $B$ nearest neighbors instead of just
one neighbor or   to 
$k$ randomly chosen sites 
\cite{gade1999synchronization,gade1996synchronization}  
the eigenvalues of the connectivity matrix can be found analytically. The stability
conditions for synchronized fixed point can be studied even in such
cases.
If the eigenvalues can be determined only numerically, we can still use the
stability conditions to explore the stability of the system by
systematically increasing the system size.

\section{Fractional order coupled map lattices: Nonlinear systems}
\begin{figure}%
	\centering
	\includegraphics[width=10cm]{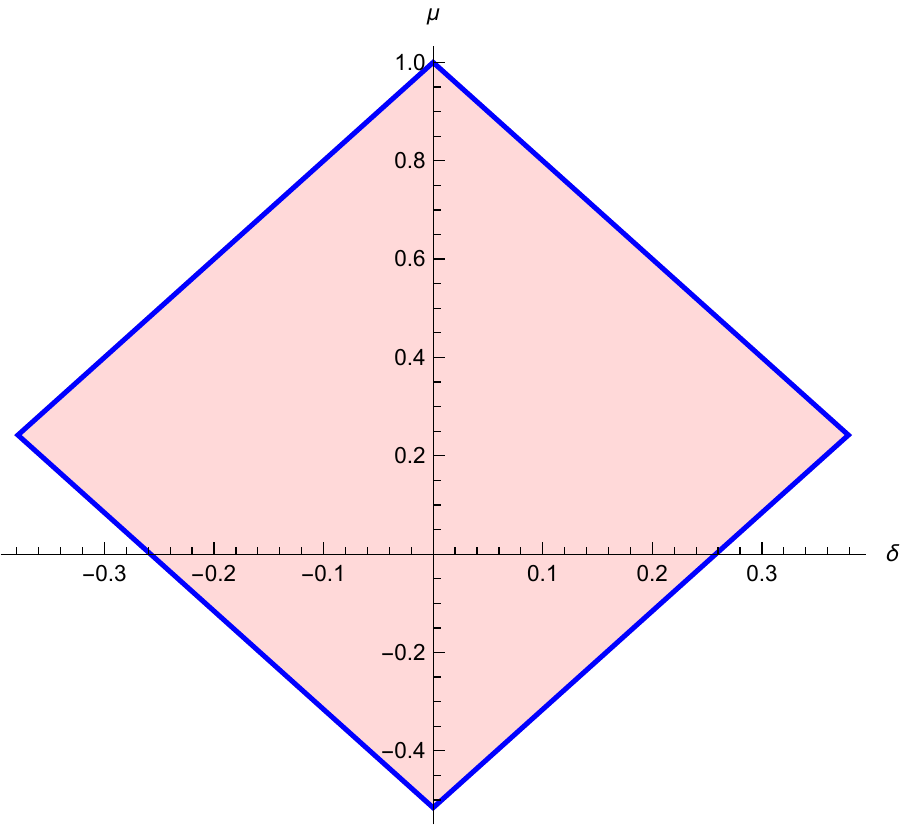} 
	\caption{Stable region of origin of (\ref{nlinCMLSys}) in Ex. \ref{nlinex1} with $N=4$ and $\alpha=0.6$}%
	\label{fig15}%
\end{figure}

\begin{figure}%
	\centering
	\includegraphics[width=10cm]{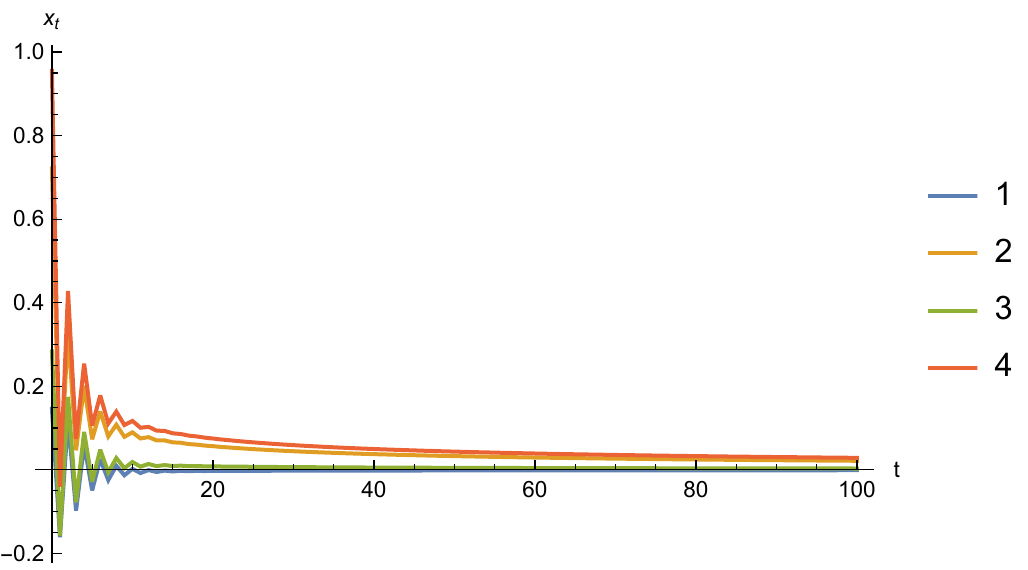} 
	\caption{Stable orbits of system (\ref{nlinCMLSys}) in Ex. \ref{nlinex1} with  $\mu=0.05$, $\delta=-0.1$, $N=4$ and $\alpha=0.6$}%
	\label{fig16}%
\end{figure}

 Consider the nonlinear coupled map lattice of fractional order $\alpha\in(0, 1]$
 \begin{equation}
 x_{t+1}(k)=x_0(k)+\sum_{j=0}^{t}\frac{\Gamma(t-j+\alpha)}{\Gamma(\alpha) \Gamma(t-j+1)}\left(f_0\left(x_j(k-1)\right)+f_1\left(x_j(k)\right)-x_j(k)+f_2\left(x_j(k+1)\right)\right), \label{nlinCML}
 \end{equation}
 where $k=1,2,\cdots,N$, $x_t(0)=x_t(N)$, $x_t(N+1)=x_t(1)$ and the functions $f_k:\mathbb{R}\longrightarrow \mathbb{R}$, $k=0,1,2$ are continuously differentiable functions.\\
If we define $X_t$ as in Section \ref{lin} and $F:\mathbb{R}^N\longrightarrow \mathbb{R}^N$ as
\begin{equation*}
F(X_t)=\begin{bmatrix}
f_0\left(x_j(N)\right)+f_1\left(x_j(1)\right)+f_2\left(x_j(2)\right)\\
f_0\left(x_j(1)\right)+f_1\left(x_j(2)\right)+f_2\left(x_j(3)\right)\\
f_0\left(x_j(2)\right)+f_1\left(x_j(3)\right)+f_2\left(x_j(4)\right)\\
\vdots\\
f_0\left(x_j(N-2)\right)+f_1\left(x_j(N-1)\right)+f_2\left(x_j(N)\right)\\
f_0\left(x_j(N-1)\right)+f_1\left(x_j(N)\right)+f_2\left(x_j(1)\right)
\end{bmatrix}
\end{equation*}
then the system (\ref{nlinCML}) is equivalent to
\begin{equation}
X_{t+1}=X_0 + \sum_{j=0}^{t}\frac{\Gamma(t-j+\alpha)}{\Gamma(\alpha) \Gamma(t-j+1)} \left[F\left(X_j\right)-X_j\right].\label{nlinCMLSys}
\end{equation}
A point $X_*=\left(x_*(1), x_*(2), \cdots, x_*(N)\right)$ is called an equilibrium point of (\ref{nlinCMLSys}) if it is a fixed point of function $F$ \cite{pakhare2022synchronization}. Therefore, such a point must satisfy
\begin{equation}
f_0\left(x_*(j-1)\right)+f_1\left(x_*(j)\right)+f_2\left(x_*(j+1)\right)=x_*(j), \quad j=1, 2, \cdots, N. \label{equilcondn}
\end{equation}
For simplicity, we assume that the equilibrium point is homogeneous, i.e. $X_*=\left(x_*, x_*, \cdots, x_*\right)$ so that the conditions (\ref{equilcondn}) get reduced to a single condition
\begin{equation}
f_0\left(x_*\right)+f_1\left(x_*\right)+f_2\left(x_*\right)=x_*. \label{hequilcondn}
\end{equation}
If we identify $a_0=f_0'\left(x_*\right)$, $a_1=f_1'\left(x_*\right)$ and $a_2=f_2'\left(x_*\right)$ then the linearization of (\ref{nlinCMLSys}) at homogeneous equilibrium point $X_*$ is given by the equation (\ref{linCMLSys}). Furthermore, if we assume the condition (\ref{hequilcondn}) then all the stability results viz. Theorems \ref{stabth0}, \ref{stabth1}, \ref{sym} and \ref{asymStabr} can be used to analyze the stability of $X_*$. We illustrate these results in the following examples.
\begin{Ex}\label{nlinex1}
Consider $f_1(x)=\mu x(1-x)$, the logistic map \cite{robert1976simple} and $f_2(x)=f_0(x)=4 x^3 - \delta x$.
\end{Ex}
Here, the origin  $X_*=(0,0,\cdots,0)$ is an equilibrium point of (\ref{nlinCMLSys}). Further, $a_1=f_1'(0)=\mu$ and $a_0=a_2=f_2'(0)=-\delta$. Therefore, for $\alpha=0.6$ and $N=4$ the stable region in $\delta\mu$-plane is sketched in Figure \ref{fig15}.
The stable orbits for $\mu=0.05$, $\delta=-0.1$ are plotted in Figure \ref{fig16}.

\begin{figure}%
	\centering
	\includegraphics[width=10cm]{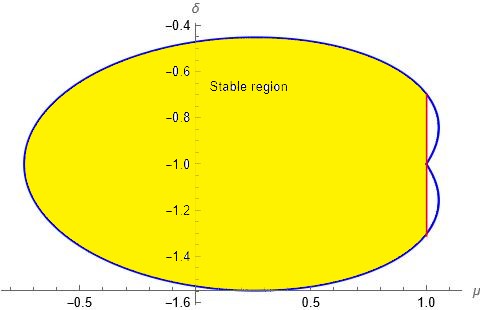} 
	\caption{Stable region of origin of (\ref{nlinCMLSys})  in Ex. \ref{nlinex2} with $N=7$ and $\alpha=0.8$}%
	\label{fig17}%
\end{figure}

\begin{figure}%
	\centering
	\includegraphics[width=10cm]{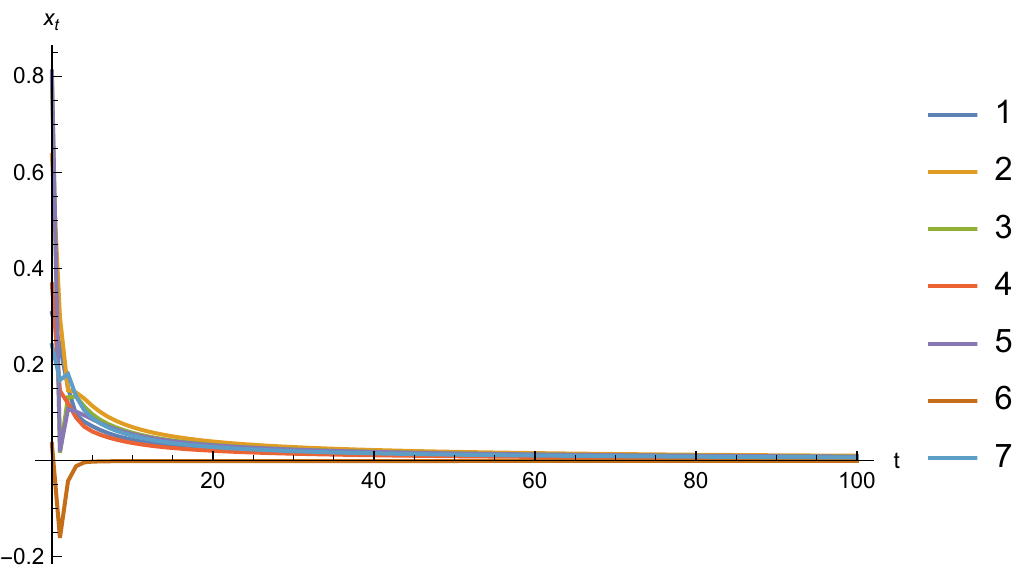} 
	\caption{Stable orbits of system (\ref{nlinCMLSys})  in Ex. \ref{nlinex2} with $\mu=0.6, \delta=-0.8$, $N=7$ and $\alpha=0.8$}%
	\label{fig18}%
\end{figure}

\begin{figure}%
	\centering
	\includegraphics[width=10cm]{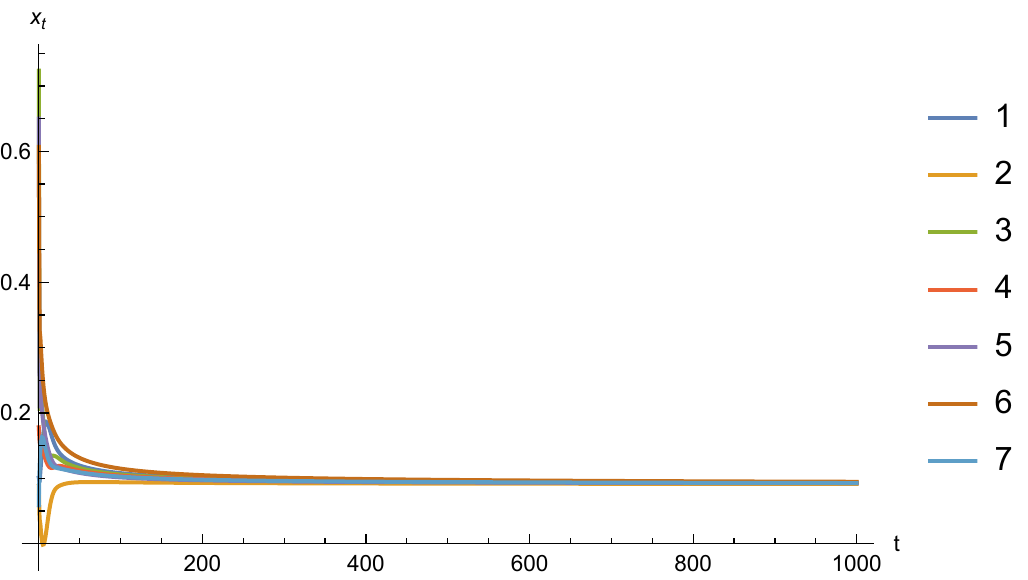} 
	\caption{Orbits of system (\ref{nlinCMLSys}) in Ex. \ref{nlinex2} diverging from origin for $\mu=1.1$, $\delta=-1.2$, $N=4$ and $\alpha=0.6$}%
	\label{fig19}%
\end{figure}

\begin{Ex}\label{nlinex2}
We consider $f_1(x)=\mu x(1-x)$ and the circle map \cite{glass1982fine} $f_2(x)=x+\delta \sin(x)$. We also set $f_0(x)=-f_2(x)$ so that the system (\ref{nlinCMLSys}) is asymmetric.
\end{Ex}
Again, we have origin as equilibrium $X_*$ and $a_1=\mu$, $a_2=-a_0=1+\delta$. We take $\alpha=0.8$ and $N=7$. The stable region shown in Figure \ref{fig17} is bounded by the line $\mu=1$ and the cardioid $\gamma_2$. The parameter values $\mu=0.6, \delta=-0.8$ in the stable region give rise to stable orbits as shown in Figure \ref{fig18}. If we take $\mu=1.1$ and $\delta=-1.2$ in the unstable region, then the trajectories repelled by origin are attracted by another homogeneous equilibrium point with $x_*=1-1/\mu$ for the sufficiently small positive initial conditions (cf. Figure \ref{fig19}). Note that the trajectories will be unbounded if we take negative initial conditions, in this case.

\begin{Ex}
Consider $f(x)=\mu x(1-x)$ and define $f_1(x)=(1-\epsilon)f(x)$ and $f_0(x)=f_2(x)=\frac{\epsilon}{2}f(x)$.
\end{Ex}
If $x_*$ is a fixed point of $f$ then $f(x_*)=x_*$ and hence the condition (\ref{hequilcondn}) is satisfied. Therefore, for this choice of functions the system (\ref{nlinCMLSys}) will have two homogeneous equilibrium points viz. $X_{1*}=(0,0,\cdots,0)$ and $X_{2*}=(q,q,\cdots,q)$, where $q=\frac{\mu-1}{\mu}$.\\
\textbf{Stability of $X_{1*}$:}\\
Here, $a_1=f_1'(0)=\mu(1-\epsilon)$ and $a_2=f_2'(0)=\epsilon\mu/2$. Therefore, $\epsilon=\frac{2a_2}{a_1+2a_2}$ and $\mu=a_1+2a_2$.
The stable region of $X_{1*}$ in the $\epsilon\mu$-plane can now be obtained using the Theorem \ref{sym} by substituting the values of $a_1$ and $a_2$ in the expressions of $\epsilon$ and $\mu$ for various values of $N$ and $\alpha$.\\
\textbf{Stability of $X_{2*}$:}\\
In this case, $a_1=f_1'(q)=(1-\epsilon)(2-\mu)$ and $a_2=f_2'(q)=\epsilon (2-\mu)/2$. On simplifying, we get  $\epsilon=\frac{2a_2}{a_1+2a_2}$ and $\mu=2-a_1-2a_2$. The stable region of $X_{2*}$ can now be traced in $\epsilon\mu$-plane by utilizing Theorem \ref{sym}.\\
The asymmetric case $f_0(x)=-f_2(x)$ can also be done in a similar way.

\section{Discussion}

As mentioned above, if the eigenvalues of underlying connectivity
matrix can be found analytically, the stability of
the synchronized state becomes very simple even
for  coupled fractional maps with an altered stability condition.
One possible extension is stability analysis of spatially 
periodic fixed point.
If an unsynchronized but spatially periodic fixed point is realized 
in fractional coupled maps (which is possible only in nonlinear systems),
the Jacobian can be 
block diagonalized. These blocks have a dimension of periodicity 
in space \cite{gade1993spatially}.  This simplifies the stability
analysis considerably.

Transition to a frozen or absorbing state is an extensively studied
transition in nonequilibrium statistical
physics which
includes systems such as coupled
oscillators. (Such transition is not possible in equilibrium systems 
because detailed balance cannot be violated.)
The above work allows us
to study such dynamical  systems in presence of memory. The
thermodynamic and asymptotic limit is important because phase can 
be defined only for the state of an infinite system  after infinite
time. The above analysis gives an analytic estimate for critical point for
such system and also gives important information about the nature
of instability. Of course, such systems can have a very different nature.
For coupled fractional maps, a power-law decay is obtained
throughout the absorbing phase and not just the critical point \cite{pakhare}.
Thus nature of transition can be very different.
In integer order maps, the bifurcation depends on whether the
eigenvalue crosses the unit circle at 1, -1, or complex value \cite{strogatz}. 
It also depends on which eigenmodes become unstable \cite{wavegade}.
Similar studies can be carried out in fractional systems.

\section*{Acknowledgment}
S. Bhalekar acknowledges the University of Hyderabad for Institute of Eminence-Professional Development Fund (IoE-PDF) by MHRD (F11/9/2019-U3(A)). P. M. Gade thanks DST-SERB for financial assistance (Ref. CRG/2020/003993).

\bibliographystyle{elsarticle-num}
\bibliography{reference.bib}
\end{document}